\documentclass[a4paper,14pt]{article}
\usepackage[english,russian]{babel}
\usepackage{amsmath}
\usepackage{amssymb}
\usepackage{amsfonts}
\usepackage{graphicx}
\date{}
\usepackage[margin=3cm]{geometry}

\begin{document}

\title {Point classification of the second order ODE's  and its application to Painleve equations}

\author{Vera V. Kartak\\
Z.Validi 32, 450074 Ufa, Russia,\\
kvera@mail.ru}
\maketitle

\noindent{\bf 2000 Mathematics Subject Classification: }{53A55, 34A26, 34A34, 34C14, 34C20, 34C41}

\noindent{\bf Key words:} Invariant, Problem of equivalence, Point transformation, Painleve equation

\maketitle

\noindent {\bf Abstract.} The first part of this work is a review of the point classification of second order ODEs done by Ruslan Sharipov. His works  were published in 1997-1998 in the Electronic Archive at LANL. The second part  is an application of this classification to Painlev\'e equations. In particular, it allows us to %
 solve  the equivalence problem for  Painlev\'e equations in an algorithmic form.

\section{Introduction}

It is a well-known fact that the following class of the second order ODEs
\begin{equation}\label{eq}
y''=P(x,y)+3\,Q(x,y)y'+3\,R(x,y)y^{\prime 2}+S(x,y)y^{\prime 3}
\end{equation}
is closed  under the generic point transformations
\begin{equation}\label{zam}
\tilde x=\tilde x(x,y),\quad \tilde y=\tilde y(x,y).
\end{equation}
It means that the transformed equation is again given by (\ref{eq}) but with some other coefficients:
\begin{equation}\label{eqt}
\tilde y''=\tilde P(\tilde x,\tilde y)+3\,\tilde Q(\tilde x,\tilde y)\tilde y'+3\,\tilde R(\tilde x, \tilde y)\tilde y^{\prime 2}+\tilde S(\tilde x,\tilde y)\tilde y^{\prime 3}.
\end{equation}
Suppose we are given two arbitrary equations (\ref{eq}) and (\ref{eqt}). The problem of
existence of the change of variables (\ref{zam}) that transforms equations (\ref{eq}) and (\ref{eqt}) one into the other  is called {\it the Equivalence Problem}. If we apply  transformation (\ref{zam}) for  equation (\ref{eq}), we  get the explicit formulas for the coefficients $\tilde P(\tilde x(x,y),\tilde y(x,y))$, $\tilde Q(\tilde x(x,y),\tilde y(x,y))$, $\tilde R(\tilde x(x,y),\tilde y(x,y))$, and $\tilde S(\tilde x(x,y),\tilde y(x,y))$ in terms of %
 $P(x,y)$, $Q(x,y)$, $R(x,y)$, $S(x,y)$ and the partial derivations of the unknown functions $\tilde x(x,y)$ and $\tilde y(x,y)$ on $x$ and $y$ up to the third order. These formulas are rather complicated, and in general situation 
the equivalence problem 
can not be solved  explicitly.

The main approach usually employed %
is to find 
invariants of equations (\ref{eq}).
{\it Invariant}  is a  function  that is preserved by transformations (\ref{zam}), i.e., %
  $I(x,y)= I(\tilde x (x,y),\tilde y(x,y)).$
Invariants 
 Theory of equations (\ref{eq}) was initiated in the  works of R.Liouville \cite {Liouville}, S.Lie \cite{Lie}, A.Tresse \cite{Tresse1, Tresse2}, E.Cartan \cite{Cartan, Thomsen}  (Late 19th- and Early 20th-Century) and continued in the Late 20th-Century in works \cite{Grissom, Lamb, Hietarinta, BordagDruima, Babich} and others. Background is  described  in papers of L.~Bordag \cite{Babich, BordagBandle}.

However, only advanced computer software for symbolic calculations gave an opportunity  to make a
 substantial progress. In the series of papers \cite{Sharipov1, Sharipov2, Sharipov3}
Ruslan Sharipov succeeded to construct the system of  (pseudo)invariants which he %
calculated explicitly in the terms of the coefficients of equations (\ref{eq}).
On the basis of this system he classified equations (\ref{eq}).
This classification is more general than all previous ones. The relation between the (pseudo)invariants from  works \cite{Sharipov2, Sharipov3} and the semiinvariants from works \cite{Cartan, Liouville} (as they were presented in  \cite{BordagBandle}) was shown in 
paper \cite{Kartak2} and here in Section 7.  Moreover, in all possible cases the set
of the invariants can be broadened.
By employing  this technique, in \cite{Kartak1}, \cite{Kartak2} and \cite{Kartak3} the equivalence problem for some equations was solved.

The first part of the present paper is a survey of  
\cite{Sharipov1, Sharipov2, Sharipov3}. We also add  additional subcases (see Subsection 5.8) not mentioned in the cited works.
The second part is an application of this classification for studying Painlev\'e equations.

\section{Classification}

  \textit{Pseudoinvariant  of weight} $m$  is a  function  transformed under  transformations (\ref{zam}) with the factor $\det T$ (the Jacobi determinant) in the power $m$,
$$J(x,y)=(\det T)^m \cdot  J(\tilde x (x,y),\tilde y(x,y)),\quad
T=\left(\begin{array}{cc}
	 \partial \tilde x /\partial  x& \partial \tilde x /\partial  y\\ \partial \tilde y /\partial  x & \partial \tilde y /\partial  y
	 \end{array}\right).
$$

\textit{Pseudotensorial field of weight $ m $ and  valence $ (r, s) $} is an indexed set  transformed  under change of variables (\ref{zam}) by the rule
$$
F^{i_1\dots i_r}_{j_1\dots j_s}=(\det T)^m{\sum_{p_1\dots p_r}}
{\sum_{q_1\dots q_s}} S^{i_1}_{p_1}\dots S^{i_r}_{p_r}T^{q_1}_{j_1}
\dots T^{q_s}_{j_s} \tilde F^{p_1\dots p_r}_{q_1\dots q_s},\quad\text{where}\quad S=T^{-1}.
$$

Given the coefficients $P$, $Q$, $R$, and $S$  of equation (\ref{eq}), we introduce a three-dimensional array by the rule
$$
\Theta_{111}=P, \qquad \Theta_{121}=\Theta_{211}=\Theta_{112}=Q,\qquad \Theta_{122}=\Theta_{212}=\Theta_{221}=R,\qquad \Theta_{222}=S.
$$
As  ``Gramian matrices'' we take the following two, 
$$
d^{ij}=d_{ij}=\left\|\begin{array}{rr}  0 & 1\\ -1 & 0
\end{array}\right\|,\quad \text{$d^{ij}$ is a pseudotensorial field of  weight 1},
\quad  \text{$d_{ij}$ is of  weight -1.}
$$
We raise the first index
\begin{equation}\label{Theta1}
{\Theta_{ij}^k}=\sum_{r=1}^2d^{kr}\Theta_{rij}.
\end{equation}
Under the change of variables  (\ref{zam}) the quantities  $\Theta_{ij}^k$ are transformed ``almost'' as an affine connection (for transformation rule see \cite{Sharipov1}).

Using ${\Theta_{ij}^k}$ as the affine connection, we construct the ``curvature tensor''
$$
{\Omega_{rij}^k}=\frac{\partial\Theta^k_{jr}}{\partial u^i}-\frac{\partial \Theta^k_{ir}}{\partial u^j}+
\sum_{q=1}^2\Theta_{iq}^k\Theta_{jr}^q-\sum_{q=1}^2\Theta^k_{jq}\Theta^q_{ir},\quad\text{here}\; u^1=x,\, u^2=y,
$$
and the ``Ricci tensor'' $\Omega_{rj}=\sum_{k=1}^2\Omega^k_{rkj}.$
Both these objects are not  tensors.
On the contrary, %
the three-dimensional array
$$
{W_{ijk}}=\nabla_i\Omega_{jk}-\nabla_j\Omega_{ik}
$$
is a tensor. Here we employ
$\Theta_{ij}^k$ in covariant differentiation instead of  the affine connection.

Using the tensor $W_{ijk}$, we introduce two extra pseudocovectorial fields,
$$
\aligned
&\alpha_k=\frac 12\sum_{i=1}^2\sum_{j=1}^2W_{ijk}d^{ij}&&
\text{is a pseudocovectorial field of weight }1,\\
&\beta_i=3\nabla_i\alpha_kd^{kr}\alpha_r+\nabla_r\alpha_kd^{kr}\alpha_i&&
\text{is a pseudocovectorial field of weight }3.
\endaligned
$$
The pseudovectorial fields are  $\alpha^j=d^{jk}\alpha_k$ of weight $2$ and
$\beta^j=d^{ji}\beta_i$ of weight $4$.

There are only three possible cases:
\begin{enumerate}
 \item \textit{maximal degeneration case}, in which  $\boldsymbol{\alpha}$=0;
 \item \textit{intermediate  degeneration case}, in which  $3F^5=\alpha^i\beta_i=0$, i.e. the fields  $\boldsymbol{\alpha}$ and  $\boldsymbol{\beta}$ are collinear;
\item \textit{general case}, in which $3F^5=\alpha^i\beta_i\ne 0$, i.e. the fields $\boldsymbol{\alpha}$ and  $\boldsymbol{\beta}$ are non-collinear.
\end{enumerate}

\section{Maximal degeneration case}
The coordinates of the pseudovectorial field ${\boldsymbol \alpha}$ are
$\alpha^1=B$, $\alpha^2=-A$,
where
\begin{equation}\label{alpha}
\aligned A=P_{ 0.2}&-2Q_{ 1.1}+R_{ 2.0}+ 2PS_{ 1.0}+SP_{
1.0}-3PR_{ 0.1}-3RP_{ 0.1} -3QR_{ 1.0} +6QQ_{ 0.1}, \\
 B=S_{ 2.0}&-2R_{ 1.1}+Q_{ 0.2}- 2SP_{0.1}-
PS_{ 0.1}+3SQ_{ 1.0}+3QS_{ 1.0}+ 3RQ_{ 0.1}-6RR_{ 1.0}.
\endaligned
\end{equation}
Hereinafter  symbol $K_{i.j}$ denotes the partial differentiation:   $K_{i.j}={\partial ^{i+j}K}/{\partial x^i\partial y^j}.$
In this case both the conditions $A=0$ and $B=0$  hold.

{\bf Theorem 1 (Lie)}. {\it All equations  (\ref{eq}) with $A=0$ and $B=0$, where $A$ and $B$ given by (\ref{alpha}), are  equivalent to $ \tilde y''=0$  by  point transformation (\ref{zam}).
They have an 8-dimension point symmetries algebra.}

 For the details see papers \cite{Liouville}, \cite{Grissom}, and others.

\section{General case}
The pseudovectorial fields $\boldsymbol{\alpha}=(B,\,-A)$ and  $\boldsymbol{\beta}=(G,\,H)$ are non-collinear, so their scalar product is non-zero. 
The pseudoinvariant $F$ defined as 
\begin{equation}\label{F}
\aligned
3F^5&=AG+BH,\qquad\text{where A and B are from (\ref{alpha}),}\\
 G&=-BB_{ 1.0}-3AB_{ 0.1}+4BA_{ 0.1}+
3SA^2-6RBA+3QB^2,\\
H&=-AA_{ 0.1}-3BA_{ 1.0}+4AB_{ 1.0}-
3PB^2+6QAB-3RA^2,
\endaligned
\end{equation}
is of weight 5. 
Since $F\ne 0$, the functions
$\varphi_1=-{\partial \ln F}/{\partial x}$ and $\varphi_2=-{\partial \ln F}/{\partial y}$
are well-defined.
Employing  $\Theta^k_{ij}$ in relations (\ref{Theta1}), we construct  an affine connection $\Gamma^k_{ij}$  and two non-collinear vectorial fields $\boldsymbol{X}$ and $\boldsymbol{Y}$ 
$$
{\Gamma^k_{ij}}=\Theta^k_{ij}-\frac{\varphi_k\delta^k_j+\varphi_k\delta^k_i}{3},\qquad
\boldsymbol{X}=\frac{\boldsymbol{\alpha}}{F^2},\qquad
\boldsymbol{Y}=\frac{\boldsymbol{\beta}}{F^4}.
$$
Their covariant derivatives are 
linear combinations of the basis fields  $\boldsymbol{X}$ and $\boldsymbol{Y}$,
$$
\begin{aligned}
&\nabla_{\boldsymbol{X}}\boldsymbol{X}=\hat\Gamma^1_{11}\boldsymbol{X}+\hat\Gamma^2_{11}\boldsymbol{Y},\qquad
\nabla_{\boldsymbol{X}}\boldsymbol{Y}=\hat\Gamma^1_{12}\boldsymbol{X}+\hat\Gamma^2_{12}\boldsymbol{Y},\\
&\nabla_{\boldsymbol{Y}}\boldsymbol{X}=\hat\Gamma^1_{21}\boldsymbol{X}+\hat\Gamma^2_{21}\boldsymbol{Y},\qquad
\nabla_{\boldsymbol{Y}}\boldsymbol{Y}=\hat\Gamma^1_{22}\boldsymbol{X}+\hat\Gamma^2_{22}\boldsymbol{Y}.
\end{aligned}
$$
Here  $\hat\Gamma^k_{ij}$ are  scalar invariants of  equation ({\ref{eq}). In  paper \cite{Sharipov2} they were denoted by
$$
I_3=\hat\Gamma^1_{12},\quad I_6=\hat\Gamma^2_{21},\quad I_7=\hat\Gamma^1_{22},\quad
I_8=\hat \Gamma^2_{22}.
$$
Differentiating these invariants along vector fields ${\boldsymbol{X}}$ and ${\boldsymbol{Y}}$ produces  
more invariants
$$
{\boldsymbol{X}}I_k=I_{k+8},\qquad {\boldsymbol{Y}}I_k=I_{k+16}.
$$
Repeating the procedure of differentiation along ${\boldsymbol{X}}$ and ${\boldsymbol{Y}}$,
we can construct an infinite sequence of invariants.
The explicit formulas for the basic four invariants read as
$$
\begin{aligned}
I_3&=\frac{B(HG_{1.0}-GH_{1.0})}{3F^9}-\frac{A(HG_{0.1}-GH_{0.1})}{3F^9}+
\frac{HF_{0.1}+GF_{1.0}}{3F^5}+\\
&+\frac{BG^2P}{3F^9}-\frac{(AG^2-2HBG)Q}{3F^9}+\frac{(BH^2-2HAG)R}{3F^9}-
\frac{AH^2S}{3F^9},\\
I_6&=\frac{H(AB_{0.1}-BA_{0.1})}{3F^7}+\frac{G(AB_{1.0}-BA_{1.0})}{4F^7}-
\frac{(AF_{0.1}-BF_{1.0})}{3F^3}-\\
&-\frac{GB^2P}{3F^7}-\frac{(HB^2-2GBA)Q}{3F^7}-\frac{(GA^2-2HBA)R}{3F^7}-
\frac{HA^2S}{3F^7},\\
I_7&=\frac{GHG_{1.0}-G^2H_{1.0}+H^2G_{0.1}-HGH_{0.1}+G^3P+3G^2HQ+3GH^2R+H^3S}{3F^{11}},\\
I_8&=\frac{G(AG_{1.0}+BH_{1.0})}{3F^9}+\frac{H(AG_{0.1}+BH_{0.1})}{3F^9}-
\frac{10(HF_{0.1}+GF_{1.0})}{3F^5}-\\
&-\frac{BG^2P}{3F^9}+\frac{(AG^2-2HBG)Q}{3F^9}-\frac{(BH^2-2HAG)R}{3F^9}+
\frac{AH^2S}{3F^9}.
\end{aligned}
$$
The case of general position splits into three subcases:
\begin{enumerate}
\item in the infinite sequence of invariants $I_k$ there exist two functionally independent ones; in this case the dimension of the point symmetries algebra is $dim(Z)=0$;
\item invariants $I_k$ are functionally dependent, but not all of them are constants; in this case $dim(Z)=1$;
\item all invarians in the sequence $I_k$ are constants; here $dim(Z)=2$.
\end{enumerate}

{\bf Example}. For equation (6.54) in the handbook by E.~Kamke \cite{Kamke}
$$
y''=y^2+4yy'+y^2y^{\prime 2}\qquad \text {we have }\quad dim(Z)=1.
$$

\section{Intermediate degeneration case}

In this case $F=0$, but $A\ne 0$ or $B\ne 0$, and the pseudovectorial fields $\boldsymbol{\alpha}$ and  $\boldsymbol{\beta}$ are collinear.

 In the case $A\ne 0$ by $\varphi_1$ and $\varphi_2$ we redenote the functions
  \begin{equation}\label{phi1}
\varphi_1=-3\frac {BP+A_{ 1.0}}{5A}+\frac 35Q, \;\;
\varphi_2=3B\frac {BP+A_{ 1.0}}{5A^2} -3\frac {B_{ 1.0}+A_{
0.1}+3BQ}{5A}+\frac 65R,
\end{equation}
and in the case $B\ne 0$ we let
\begin{equation}\label{phi2}
 \varphi_1=-3A\frac {AS-B_{ 0.1}}{5B^2}
-3\frac {A_{ 0.1}+B_{ 1.0}-3AR}{5B}-\frac 65Q,\quad
\varphi_2=3\frac {AS-B_{ 0.1}}{5B}-\frac 35R.
\end{equation}
Employing the introduced functions and $\Theta^k_{ij}$ from (\ref{Theta1}), we
construct the affine connection $\Gamma^k_{ij}$    and a pseudoinvariant $\Omega$ of weight 1,
\begin{equation}\label{conn}
{\Gamma^k_{ij}}=\Theta^k_{ij}-\frac{\varphi_k\delta^k_j+\varphi_k\delta^k_i}{3},\qquad
{\Omega}=\frac 53\left( \frac{\partial\varphi_1}{\partial y}-\frac{\partial\varphi_2}{\partial x}
\right).
\end{equation}
As $A\ne 0$, the explicit formula for the pseudoinvariant  $\Omega$ reads as
\begin{equation}\label{Omega1}
\aligned \Omega &=\frac {2BA_{ 1.0}(BP+ A_{ 1.0})}{A^3}- \frac
{(2B_{ 1.0}+3BQ)A_{ 1.0}}{A^2}+\frac {(A_{ 0.1}-2B_{1.0})BP}{A^2}-\\
&- \frac {BA_{ 2.0}+B^2 P_{ 1.0}}{A^2}+ \frac
{B_{ 2.0}}A+\frac {3B_{ 1.0}Q+3BQ_{ 1.0}- B_{ 0.1}P-BP_{
0.1}}{A}+Q_{ 0.1}- 2R_{ 1.0}.
\endaligned
\end{equation}
And in the case  $B\ne 0$ the similar formula is
\begin{equation}\label{Omega2}
\aligned \Omega &=\frac {2AB_{ 0.1}(AS- B_{ 0.1})}{B^3}- \frac
{(2A_{ 0.1}-3AR)B_{ 0.1}}{B^2}+\frac {(B_{ 1.0}-2A_{
0.1})AS}{B^2}+ \\
&+\frac {AB_{ 0.2}-A^2 S_{ 0.1}}{B^2}- \frac
{A_{ 0.2}}B+\frac {3A_{ 0.1}R+3AR_{ 0.1}- A_{ 1.0}S-AS_{
1.0}}{B}+R_{ 1.0}- 2Q_{ 0.1}.
\endaligned
\end{equation}

The rule of covariant differentiation of the pseudotensorial field was given in \cite{Sharipov1},
$$
\nabla_kF^{i_1\dots i_r}_{j_1\dots j_s}= \frac{\partial F^{i_1\dots i_r}_{j_1\dots j_s}}{\partial u^k}+
\sum_{n=1}^r\sum_{v_n=1}^2\Gamma_{k v_n}^{i_n}F^{i_1\dots v_n\dots i_r}_{j_1\dots j_s}
-\sum_{n=1}^s\sum_{w_n=1}^2\Gamma_{k j_n}^{w_n}F^{i_1\dots i_r}_{j_1\dots w_n\dots  j_s}+
m\varphi_k F^{i_1\dots i_r}_{j_1\dots j_s}.
$$
If a pseudotensorial field $F$ has valence $(r,s)$ and weight $m$, the
pseudotensorial field $\nabla F$ has valence $(r, s+1)$ and weight $m.$

The pseudovectorial fields  $\boldsymbol{\alpha}$ and  $\boldsymbol{\beta}$ are collinear, hence there exists a coefficient $N$ such that
$\boldsymbol{\beta}=3N\boldsymbol{\alpha}.$
$N$  is a pseudoinvariant of  weight 2.

We let
\begin{equation}\label{xi}
{\xi^i}=d^{ij}\nabla_jN,\qquad {M}=-\alpha_i\xi^i,\qquad \boldsymbol{\gamma}=-\boldsymbol{\xi}-2\Omega\boldsymbol{\alpha},
\end{equation}
Here  $\boldsymbol{\xi}$ is a pseudovectorial field of  weight 3,  ${M}$ is a pseudoinvariant of weight 4,
$\boldsymbol{\gamma}$ is a pseudovectorial field of  weight 3.

In the cases $A\ne 0$ and $B\ne 0$ the pseudoinvariant $N$ is given by the formulas
\begin{equation}\label{N}
 N =-\frac H{3A}, \qquad\qquad N =\frac G{3B},
\end{equation}
respectively. 
The pseudoinvariant $M$  in the case $A\ne 0$ reads as
\begin{equation}\label{M1}
M=-\frac {12BN(BP+A_{ 1.0})}{5A}+BN_{ 1.0}+\frac {24}5BNQ+\frac
65NB_{ 1.0}+\frac 65NA_{ 0.1}-AN_{ 0.1}- \frac {12}5ANR,
\end{equation}
and in the case $B\ne 0$ it is given by the formula
\begin{equation}\label{M2}
M=-\frac {12AN(AS-B_{ 0.1})}{5B}-AN_{ 0.1}+\frac {24}5ANR-
 \frac
65NA_{ 0.1}-\frac 65NB_{ 1.0}+BN_{ 1.0}- \frac {12}5 BNQ.
\end{equation}
In the case $A\ne 0$ the field ${\boldsymbol \gamma}$ is 
\begin{equation}\label{gamma1}
\aligned
 \gamma^1=&-\frac {6BN(BP+A_{ 1.0})}{5A^2}+
\frac {18NBQ}{5A}+
\frac {6N(B_{ 1.0}+A_{ 0.1})}{5A} -N_{ 0.1}-\frac
{12}5NR-2\Omega B,\\
\gamma^2=&-\frac {6N(BP+A_{ 1.0})}{5A}+N_{ 1.0}+\frac 65NQ+ 2\Omega
A.
\endaligned
\end{equation}
In the case $B\ne 0$ the field ${\boldsymbol \gamma}$ is
\begin{equation}\label{gamma2}
\aligned
\gamma^1=&-\frac {6N(AN-B_{ 0.1})}{5B}-N_{ 0.1} +\frac 65NR-2\Omega
B,\\
 \gamma^2=&-\frac {6AN(AS-B_{ 0.1})}{5B^2}+
\frac {18NAR}{5B}-
\frac {6N(A_{ 0.1}+B_{ 1.0})}{5B} +N_{ 1.0}-\frac
{12}5NQ+2\Omega A.
\endaligned
\end{equation}

\subsection{First case of intermediate degeneration: $M\ne 0$}
If in (\ref{M1}), (\ref{M2}) $M\ne 0$, the pseudovectorial fields ${\boldsymbol \alpha}$ in (\ref{alpha}) and ${\boldsymbol \gamma}$ in (\ref{gamma1}), (\ref{gamma2})
are non-collinear. Moreover, it means that $N\ne 0$ in (\ref{N}). Consider the  expansion
$\nabla_{\boldsymbol \gamma}{\boldsymbol \gamma}=\hat\Gamma^1_{22}{\boldsymbol \alpha}+
\hat\Gamma^2_{22}{\boldsymbol \gamma}.$
The basic invariants are the following ones,
\begin{equation}\label{inv1}
I_1=\frac{M}{N^2},\qquad I_2=\frac{\Omega^2}{N},\qquad I_3=\frac{\hat\Gamma^1_{22}}{M}.
\end{equation}
Here$M$, $N$, and $\Omega$ are from (\ref{M1}), (\ref{M2}), (\ref{N}), (\ref{Omega1}), (\ref{Omega2}). The explicit formula for $\hat\Gamma^1_{22}$ is
$$
\hat\Gamma^1_{22}=\frac {\gamma^1\gamma^2(\gamma^1_{
1.0}- \gamma^2_{ 0.1})}{M}+ \frac {(\gamma^2)^2\gamma^1_{ 0.1}-
(\gamma^1)^2\gamma^2_{ 1.0}}M+
\frac
{P(\gamma^1)^3+3Q(\gamma^1)^2\gamma^2+3R\gamma^1(\gamma^2)^2+
S(\gamma^2)^3}M.
$$
By differentiating the invariants $I_1$, $I_2$ and $I_3$ along fields ${\boldsymbol \alpha}$ and ${\boldsymbol \gamma}$ we get new invariants
\begin{equation}\label{inv1-1}
I_{k+3}=\frac{\nabla_{\boldsymbol \alpha}I_k}{N},\qquad I_{k+6}=\frac{(\nabla_{\boldsymbol \gamma}I_k)^2}{N^3}.
\end{equation}
The first case of intermediate degeneration splits into three subcases
\begin{enumerate}
\item in the infinite sequence of invariants $I_k$ there exist two functionally independent ones; in this case the dimension of the point symmetries algebra is $dim(Z)=0$;
\item invariants $I_k$ are functionally dependent but not all of them are constants; here we have $dim(Z)=1$;
\item all invariants in the sequence $I_k$ are constants; here $dim(Z)=2$.
\end{enumerate}

{\bf Example}. For equation (6.45) in the handbook by E.~Kamke \cite{Kamke}
$$
y''=ay^{\prime 2}+by,\quad ab\ne 0,\qquad \text{one has}\quad dim(Z)=1.
$$

\subsection{Second case of intermediate degeneration}
If in (\ref{M1}), (\ref{M2}) $M=0$, the pseudovectorial fields ${\boldsymbol \alpha}$ in (\ref{alpha}) and ${\boldsymbol \gamma}$ in (\ref{gamma1}), (\ref{gamma2}) are collinear. Hence, there exists a coefficient $\Lambda$  such that
$\boldsymbol{\gamma}=\Lambda\boldsymbol{\alpha}.$
Here $\Lambda$ is a pseudoinvariant of weight 1, in the cases $A\ne 0$ and $B\ne 0$ being respectively
$$
\Lambda=-\frac{\gamma^2}{A},\quad A\ne 0,\qquad\mbox{or} \qquad \Lambda=\frac{\gamma^1}{B},\quad B\ne 0.
$$
The explicit formulas for $\Lambda$ are
\begin{equation}\label{lambda1}
\Lambda=-\frac{6N(AS-B_{0.1})}{5B^2}-\frac{N_{0.1}}B+\frac{6NR}{5B}-2\Omega.
\end{equation}
\begin{equation}\label{lambda2}
\Lambda=\frac{6N(BP+B_{1.0})}{5A^2}-\frac{N_{1.0}}A-\frac{6NQ}{5A}-2\Omega.
\end{equation}

Let us calculate the curvature tensor using the connections (\ref{conn}):
$$
R^k_{qij}=\frac{\partial \Gamma^k_{jk}}{\partial u^i}-
\frac{\partial\Gamma^k_{iq}}{\partial u^j}+\sum_{s=1}^2\Gamma^k_{is}\Gamma^s_{jq}-
\sum_{s=1}^2\Gamma^k_{js}\Gamma^s_{iq},\quad u^1=x,\, u^2=y,
$$
and the pseudotensorial field of the weight 1:
$$
R^k_q=\frac 12\sum_{i=1}^2\sum_{j=1}^2 R^i_{qij}d^{ij},
$$
where $\lambda_1$ and $\lambda_2$ are its eigenvalues.
Now we construct the pseudocovectorial field of the weight -1.
If $A\ne 0$, we let
\begin{equation}\label{omega1}
\omega_1=-\frac{R^2_1}{A},\qquad \omega_2=\frac{\lambda_2-R^2_2}{A},
\end{equation}
where
$$
\begin{aligned}
\omega_1&=\frac{12PR}{5A}-\frac{54}{25}\frac{Q^2}{A}-\frac{P_{0.1}}{A}+\frac{6Q_{1.0}}{5A}-
\frac{PA_{0.1}+BP_{1.0}+A_{2.0}}{5A^2}-\\
&-\frac{2B_{1.0}P}{5A^2}+\frac{3QA_{1.0}-12PBQ}{25A^2}+
\frac{6B^2P^2+12A_{1.0}BP+6A_{1.0}^2}{25A^3},\\
\omega_2&=\frac{6\Lambda+3\Omega}{5A}+\frac{-5BP_{0.1}+6BQ_{1.0}+12RBP}{5A^2}-\frac{54}{25}\frac{BQ^2}{A^2}-\frac{12B^2PQ}{25A^3}+\frac{3BQA_{1.0}}{25A^3}-\\
&-\frac{2BB_{1.0}P+BA_{0.1}P+B^2P_{1.0}+BA_{2.0}}{5A^3}+\frac{6BA_{1.0}^2+6B^3P^2+12B^2A_{1.0}P}{25A^4}.
\end{aligned}
$$
And if $B\ne 0$,
\begin{equation}\label{omega2}
\omega_1=\frac{R^1_1-\lambda_2}B,\qquad \omega_2=\frac{R^1_2}B,\qquad\text{and}
\end{equation}
$$
\begin{aligned}
\omega_1&=-\frac{6\Lambda+3\Omega}{5B}+\frac{5AS_{1.0}-6AR_{0.1}+12QAS}{5B^2}-
\frac{54}{25}\frac{AR^2}{B^2}-\frac{12A^2SR}{25B^3}+\frac{3ARB_{0.1}}{25B^3}+\\
&+\frac{2AA_{0.1}S+AB_{1.0}S+A^2S_{0.1}-AB_{0.2}}{5B^3}
+\frac{6AB_{0.1}^2+6A^3S^2-12A^2B_{0.1}S}{25B^4},\\
\omega_2&=\frac{12SQ}{5B}-\frac{54}{25}\frac{R^2}{B}+\frac{S_{1.0}}{B}-\frac{6R_{0.1}}{5B}+
\frac{SB_{1.0}+AS_{0.1}-B_{0.2}}{5B^2}+\\
&+\frac{2A_{0.1}S}{5B^2}-\frac{3RB_{0.1}+12SAR}{25B^2}+
\frac{6A^2S^2-12B_{0.1}AS+6B_{0.1}^2}{25B^3}.
\end{aligned}
$$
We introduce one more pseudocovectorial field of weight 1,
$$
{\boldsymbol w}=N{\boldsymbol \omega}+\nabla\Lambda+\frac 13\nabla\Omega.
$$
It is collinear to the pseudovectorial field ${\boldsymbol \alpha}$ in (\ref{alpha}), and thus there exists $K$ such that
${\boldsymbol w}=K{\boldsymbol \alpha}$,
\begin{equation}\label{K1}
K=\frac{\Lambda_{1.0}+\Lambda\varphi_1}{A}+\frac{\Omega_{1.0}+\Omega\varphi_1}{3A}+
\frac{N\omega_1}{A},\qquad A\ne 0.
\end{equation}
\begin{equation}\label{K2}
K=\frac{\Lambda_{0.1}+\Lambda\varphi_2}{B}+\frac{\Omega_{0.1}+\Omega\varphi_2}{3B}+
\frac{N\omega_2}{B},\qquad B\ne 0.
\end{equation}

By ${\boldsymbol \varepsilon}$ we denote a
 pseudocovectorial field  of  weight 1,
$$
{\boldsymbol \varepsilon}=N{\boldsymbol \omega}+\nabla\Lambda.
$$
Raising indices by  matrix $d^{ij}$, we get the pseudovectorial field
${\boldsymbol \varepsilon}$ of weight 2,
\begin{equation}\label{epsilon}
\varepsilon^1=N\omega_2+\Lambda_{0.1}+\varphi_2\Lambda,\qquad
\varepsilon^2=-N\omega_1-\Lambda_{1.0}+\varphi_1\Lambda.
\end{equation}
The fields ${\boldsymbol \varepsilon}$ in (\ref
{epsilon}) and ${\boldsymbol \alpha}$ in (\ref{alpha}) are non-collinear, and we can write
$$
\nabla_{\boldsymbol \varepsilon}{\boldsymbol \varepsilon}=\hat\Gamma^1_{22}{\boldsymbol \alpha}+
\hat\Gamma^2_{22}{\boldsymbol \varepsilon}.
$$
$$
\begin{aligned}
\hat\Gamma^1_{22}&=\frac{5\varepsilon^1\varepsilon^2(\varepsilon_{1.0}^1-\varepsilon^2_{0.1})}{3N\Omega}+\frac{5(\varepsilon^2)^2\varepsilon^1_{0.1}-5(\varepsilon^1)^2\varepsilon^2_{1.0}}{3N\Omega}+\\
&+\frac{5P(\varepsilon^1)^3+15Q(\varepsilon^1)^2\varepsilon^2+15R\varepsilon^1(\varepsilon^2)^2+
5S(\varepsilon^2)^3}{3N\Omega}.
\end{aligned}
$$
We introduce  pseudoscalar fields
\begin{equation}\label{L1}
L=KN+\frac 59N+3\Lambda\Omega+\frac 79\Omega^2+2\Lambda^2.
\end{equation}
\begin{equation}\label{E1}
\begin{aligned}
E=&\hat\Gamma^1_{22}-\frac{\nabla_{\boldsymbol \varepsilon}L}{N}+\frac{4\Lambda\nabla_{\boldsymbol \varepsilon}\Lambda}{N}+\frac{17\Omega\nabla_{\boldsymbol \varepsilon}\Lambda}{6N}+
\frac{12L^2}{5N}-\frac{53L\Lambda\Omega}{5N}-\frac{48L\Lambda^2}{5N}-\\
{}&-\frac{62L\Omega^2}{15N}-\frac{8L}{3}+\frac{48\Lambda^4}{5N}+\frac{106\Lambda^3\Omega}{5N}+\frac{16\Lambda^2}{3}+\\
{}&+\frac{1163\Lambda^2\Omega^2}{60N}+\frac{137\Lambda\Omega^3}{18N}+
\frac{50\Lambda\Omega}{9}+\frac{203\Omega^2}{108}+\frac{77\Omega^4}{135N}+\frac{20N}{27}.
\end{aligned}
\end{equation}

Employing the above objects, we can define
invariants %
$$
I_1=\frac{\Lambda^{12}}{\Omega^8N^2},\qquad I_2=\frac{L^4}{N^2\Omega^4},\qquad I_3=\frac{E^6N^4}{\Omega^{20}}.
$$
Here $\Lambda$ is from (\ref{lambda1}), (\ref{lambda2}), $\Omega$ is from (\ref{Omega1}), (\ref{Omega2}), N is from (\ref{N}), L is from (\ref{L1}), E is from (\ref{E1}):
In the second case of intermediate degeneration the algebra of the point symmetries of  equation
(\ref{eq})
is 1-dimensional if and only if all invariants $I_1,\, I_2,\, I_3$ are identically constant. In other cases it is trivial.

\subsection{Third case of intermediate degeneration}
In this case $N\ne 0$ in (\ref{N}), $M=0$ in (\ref{M1}), (\ref{M2}), $\Omega=0$ in (\ref{Omega1}), (\ref{Omega2}), $\Lambda\ne 0$ in (\ref{lambda1}), (\ref{lambda2}).  Consider again the pseudocovectorial field ${\boldsymbol \omega}$ of the weight -1 from (\ref{omega1}), (\ref{omega2}). Raising indices by the matrix $d^{ij}$, we get the vector field ${\boldsymbol \omega}$,
$\omega^1=\omega_2$, $\omega^2=-\omega_1$. Since $\Lambda\ne 0$,  ${\boldsymbol \omega}$, and ${\boldsymbol \alpha}$ are non-collinear, we obtain the following relation
$$
\nabla_{\boldsymbol \omega}{\boldsymbol \omega}=\hat\Gamma^1_{22}{\boldsymbol \alpha}+
\hat\Gamma^2_{22}{\boldsymbol \omega},
$$
$$
\begin{aligned}
\hat\Gamma^1_{22}&=\frac{5\omega^1\omega^2(\omega_{1.0}^1-\omega^2_{0.1})}{\Lambda}+\frac{5(\omega^2)^2\omega^1_{0.1}-5(\omega^1)^2\omega^2_{1.0}}{\Lambda}+\\
&+\frac{5P(\omega^1)^3+15Q(\omega^1)^2\omega^2+15R\omega^1(\omega^2)^2+
5S(\omega^2)^3}{6\Lambda}.
\end{aligned}
$$
In this case we define new $L$ and $E$,
\begin{equation}\label{L2}
\begin{aligned}
& L=K+\frac 59+\frac{2\Lambda^2}{N},\quad\text{with $K$ from (\ref{K1}), (\ref{K2})}\\
& E=\hat\Gamma^1_{22}-\frac{\nabla_{\boldsymbol \omega}L}{N}+
\frac{9L^2}{5N}-\frac{2L}{N}-\frac{12L\Lambda^2}{5N^2}+\frac{7\Lambda^2}{3N^2}+
\frac{5}{9N}+\frac{63\Lambda^4}{20N^3}.
\end{aligned}
\end{equation}
Let us construct the invariants:
$$
I_1=\frac{L^8N^6}{\Lambda^{12}},\qquad I_2=\frac{EN^3}{\Lambda^4}.
$$
Here $L$, $E$ are from (\ref{L2}), $ N$ is from (\ref{N}), $\Lambda$ is from (\ref{lambda1}), (\ref{lambda2}).

In the third case of intermediate degeneration the algebra of the point symmetries of  equation (\ref{eq})
is 1-dimensional if and only if both the invariants $I_1,\, I_2$ are identically constant. In other cases it is trivial.

{\bf Example.} For Emden-Fowler equation (6.11)  with $n=-3$ in the handbook by E.~Kamke  \cite{Kamke}
$$
y''=-\frac{ax^m}{y^3},\quad a\ne 0,\quad\text{one has}\quad dim(Z)=1.
$$

\subsection{Fouth case of intermediate degeneration}
In this case $N\ne 0$ in (\ref{N}), $M=0$ in (\ref{M1}), (\ref{M2}), $\Omega=0$ in (\ref{Omega1}), (\ref{Omega2}), $\Lambda= 0$ in (\ref{lambda1}), (\ref{lambda2}), $K\ne -5/9$ in (\ref{K1}), (\ref{K2}).
Consider again the vectorial field ${\boldsymbol \omega}$, $\omega^1=\omega_2$, $\omega^2=-\omega_1$, from (\ref{omega1}), (\ref{omega2}).  Since $\Lambda= 0$,  ${\boldsymbol \omega}$, and ${\boldsymbol \alpha}$ are collinear, we can define a new scalar field $\Theta$ by the relationship ${\boldsymbol \omega}=\Theta{\boldsymbol \alpha}$,
\begin{equation}\label{theta}
\Theta=\frac{\omega_1}{A},\quad A\ne 0,\qquad \Theta=\frac{\omega_2}{B},\quad B\ne 0.
\end{equation}
The covariant differential ${\boldsymbol \theta}=\nabla\Theta$ is a pseudocovectorial field of  weight -2,
\begin{equation}\label{thetaf}
\theta_1=\Theta_{1.0}-2\varphi_1\Theta,\qquad \theta_2=\Theta_{0.1}-2\varphi_2\Theta.
\end{equation}
The corresponding pseudovectorial field of weight -1 is $\theta^1=\theta_2$, $\theta^2=-\theta_1$. Let us calculate its convolution with ${\boldsymbol \alpha}$ from (\ref{alpha}),
\begin{equation}\label{L3}
L=-\frac 59\sum_{i=1}^2\alpha_i\theta^i.
\end{equation}
And the  relation $L=K+5/9$ holds true,
where $K$ is from (\ref{K1}), (\ref{K2}).

Since $L\ne 0$,  fields ${\boldsymbol \theta}$ (\ref{thetaf}) and ${\boldsymbol \alpha}$ (\ref{alpha}) are non-collinear,
$$
\nabla_{\boldsymbol \theta}{\boldsymbol \theta}=\hat\Gamma^1_{22}{\boldsymbol \alpha}+
\hat\Gamma^2_{22}{\boldsymbol \theta},
$$
$$
\begin{aligned}
\hat\Gamma^1_{22}&=-\frac{5\theta^1\theta^2(\theta_{1.0}^1-\theta^2_{0.1})}{9L}-\frac{5(\theta^2)^2\theta^1_{0.1}-5(\theta^1)^2\theta^2_{1.0}}{9L}-\\
&-\frac{5P(\theta^1)^3+15Q(\theta^1)^2\theta^2+15R\theta^1(\theta^2)^2+
5S(\theta^2)^3}{9L}.
\end{aligned}
$$
We introduce one more pseudoscalar field,
\begin{equation}\label{E3}
E=\hat\Gamma^1_{22}+\frac{27N}{5}\left(\Theta+\frac{5}{9N}\right)^3-
\frac34\left(\Theta+\frac{5}{9N}\right)^2
\end{equation}
that gives rise to the invariant
$$
I_1=\frac{E^6N^{12}}{L^{20}}
$$
where $E$ is from (\ref{E3}), $N$ is from (\ref{N}), $L$ is from (\ref{L3}).

In the fouth case of intermediate degeneration the algebra of the point symmetries of  equation (\ref{eq})
is 1-dimensional if and only if the invariant $I_1$ is identically constant. Otherwise it is trivial.

\subsection{Fifth case of intermediate degeneration}
In this case $N\ne 0$ in (\ref{N}), $M=0$ in (\ref{M1}), (\ref{M2}), $\Omega=0$ in (\ref{Omega1}), (\ref{Omega2}), $\Lambda= 0$ in (\ref{lambda1}), (\ref{lambda2}), $K= -5/9$ in (\ref{K1}), (\ref{K2}).
All equations (\ref{eq}) are equivalent to
$$
y''=\frac 1{y^3},\quad\mbox{or another form}\quad y''=-\frac 5{4x}y'+\frac 43x^2y^{\prime 3}.
$$
The algebra of point symmetries is 3-dimensional, see also \cite{Rom}.

\subsection{Sixth case of intermediate degeneration}
In this case $N= 0$ in (\ref{N}), $\Omega\ne 0$ in (\ref{Omega1}), (\ref{Omega2}).
The pseudovectorial fields ${\boldsymbol \omega}$, $\omega^1=\omega_2$, $\omega^2=-\omega_1$ from (\ref{omega1}), (\ref{omega2}) and ${\boldsymbol \alpha}$ from (\ref{alpha}) are non-collinear,
$$
\nabla_{\boldsymbol \omega}{\boldsymbol \omega}=\hat\Gamma^1_{22}{\boldsymbol \alpha}+
\hat\Gamma^2_{22}{\boldsymbol \omega},
$$
$$
\begin{aligned}
\hat\Gamma^1_{22}&=-\frac{5\omega^1\omega^2(\omega_{1.0}^1-\omega^2_{0.1})}{9\Omega}-\frac{5(\omega^2)^2\omega^1_{0.1}-5(\omega^1)^2\omega^2_{1.0}}{9\Omega}-\\
&-\frac{5P(\omega^1)^3+15Q(\omega^1)^2\omega^2+15R\omega^1(\omega^2)^2+
5S(\omega)^3}{9\Omega}.
\end{aligned}
$$
The corresponding invariants are
$$
\begin{aligned}
&I_1=L=\nabla_{\boldsymbol \omega}K-\frac{21}{25}K^2-K,\\
&I_2=\Omega^2\hat\Gamma^1_{22}-\nabla_{\boldsymbol \omega}L-\frac{72}{625}K^3+\frac{63}{50}K^2+
\frac{12}{25}KL-K-L.
\end{aligned}
$$
where $K$ is from (\ref{K1}), (\ref{K2}), $\Omega$ is from (\ref{Omega1}), (\ref{Omega2}).

In the sixth case of intermediate degeneration the algebra of point symmetries  of  equation
(\ref{eq}) is 1-dimensional if and only if both invariants $I_1,\, I_2$ are identically constant. In other cases it is trivial.

\subsection{Seventh case of intermediate degeneration}
In this case $N= 0$ in (\ref{N}), $\Omega= 0$ in (\ref{Omega1}), (\ref{Omega2}).
The pseudovectorial fields ${\boldsymbol \theta}$ from (\ref{thetaf}) and ${\boldsymbol \alpha}$ from (\ref{alpha}) are non-collinear,
$$
\nabla_{\boldsymbol \theta}{\boldsymbol \theta}=\hat\Gamma^1_{22}{\boldsymbol \alpha}+
\hat\Gamma^2_{22}{\boldsymbol \theta},
$$
$$
\begin{aligned}
\hat\Gamma^1_{22}&=\theta^1\theta^2(\theta_{1.0}^1-\theta^2_{0.1})-(\theta^2)^2\theta^1_{0.1}+(\theta^1)^2\theta^2_{1.0}-\\
&-P(\theta^1)^3-3Q(\theta^1)^2\theta^2-3R\theta^1(\theta^2)^2-S(\theta)^3.
\end{aligned}
$$
We define a pseudoscalar field and an invariant,
\begin{equation}\label{inv7}
L=\hat\Gamma^1_{22}-\frac 12\Theta^2, \qquad I_1=\frac{L_1^4}{L^5},  \qquad I_2=\frac{\Theta^2}{L},
\end{equation}
where $\Theta$ is from (\ref{theta}), and a pseudoinvariant
\begin{equation}\label{L1}
L_1= \nabla_{\boldsymbol \theta}L=L_{{1.0}}\theta^1+ L_{
0.1}\theta^2-4L(\varphi_1\theta^1+\varphi_2\theta^2)
\end{equation}
where we have employed  (\ref{inv7}), (\ref{theta}) and (\ref{phi1}), (\ref{phi2}).

In the seventh case of intermediate degeneration the algebra of the point symmetries  of  equation (\ref{eq})
is 2-dimensional if and only if  $L=0$; and
is 1-dimensional if and only if  $L\ne 0$ and if  $I_1$ is identically constant. In other cases the algebra is trivial.

{\bf Example.} For equation (6.5) in the handbook by  E.~Kamke  \cite{Kamke}
$$
y''=ay^2+bx+c,\quad a\ne 0\quad\text{one has} \quad dim(Z)=2 \quad\text{if}\quad  b^2=ac\quad \text{otherwise} \quad dim(Z)=1.
$$

\subsection{Additional subcases of intermediate degeneration}
We define a pseudovectorial field ${\boldsymbol \eta}$, where $\eta^i=d^{ij}\nabla_jM$  and  its scalar product with the field ${\boldsymbol \xi}$ from (\ref{xi}):
\begin{equation}\label{z}
Z=d_{ij}\eta^i\xi^j.
\end{equation}
Here $Z$ is a pseudoinvariant of  weight 7.
Then the first case of the intermediate degeneration splits into  four subcases.

{\bf Subcase 1.1.} $M\ne 0$, $\Omega \ne 0$, $Z\ne 0$.

{\bf Subcase 1.2.} $M\ne 0$, $\Omega \ne 0$, $Z= 0$.

{\bf Subcase 1.3.} $M\ne 0$, $\Omega =0$, $Z\ne 0$.

{\bf Subcase 1.4.} $M\ne 0$, $\Omega = 0$, $Z= 0$.

Subject to  the pseudoinvariant $\Theta$ from (\ref{theta}), the seventh case of the intermediate degeneration  splits into the two subcases.

{\bf Subcase 7.1.} $N=0$, $\Omega =0$, $\Theta \ne 0$.

{\bf Subcase 7.2.} $N=0$, $\Omega = 0$, $\Theta= 0$.

\subsection{Tree of  intermediate degeneration cases}
The following diagram illustrates the cases of the intermediate degeneration.
\begin{figure}[h]
\center{\includegraphics[width=1\linewidth]{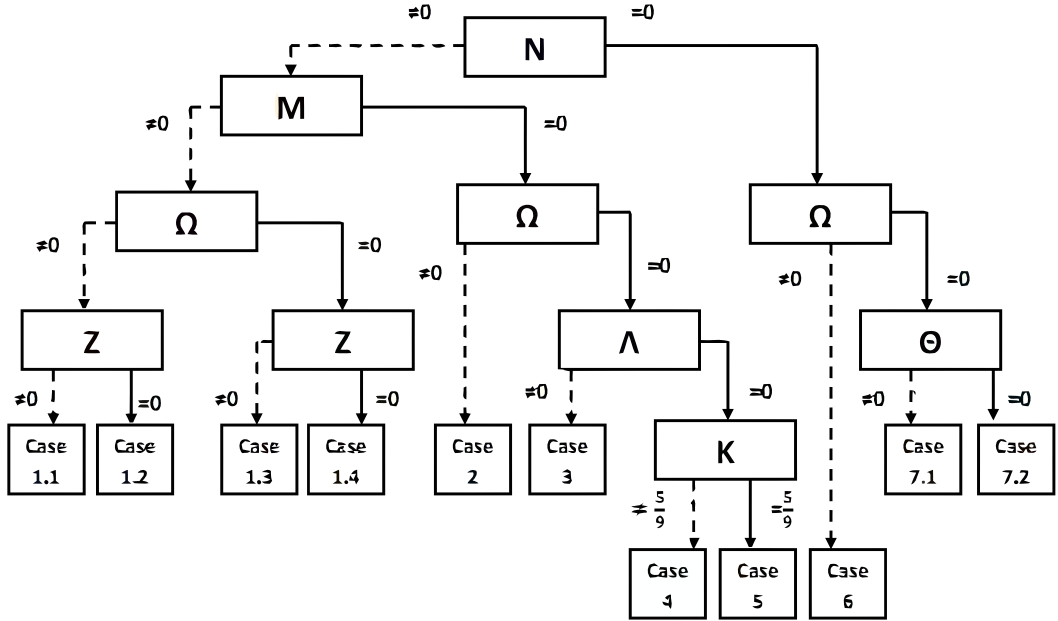}}
\vspace*{8pt}
\caption{Tree of  intermediate degeneration cases.}
\label{tree}
\end{figure}

\section{Classification of Painlev\'e equations}
Let us show how
Painlev\'e equations are included into the proposed classification scheme.

\begin{enumerate}
\item Equation Painlev\'e I  is in Case 7.1 of intermediate degeneration.

\item Equations Painlev\'e III-VI (except special cases!) are in  Case 1.3 of intermediate degeneration.

\item Special cases.

\begin{enumerate}
\item Equation Painlev\'e II is in  Case 1.4 of intermediate degenera\-tion.

\item Equation Painlev\'e III with 3 zero parameters is in Case 1.4 of intermediate degenera\-tion.

\item Equation Painlev\'e III  $(0,b,0,d)$ or $(a,0,c,0)$ (they are equivalent) is in  Case~1.4 of intermediate degeneration.

\item Equation Painlev\'e V  $(a,b,0,0)$ is in  Case 1.4 of intermediate degeneration.

\item Equation Painlev\'e III  $(0,0,0,0)$ is in  Case of maximal degeneration.

\item Equation Painlev\'e V  $(0,0,0,0)$ is in  Case of maximal degeneration.

\item Equation Painlev\'e VI  $(0,0,0,\frac 12)$ is in  Case of maximal degeneration.
\end{enumerate}
\end{enumerate}

\section{Relation between the semiinvariants}

In work of  E.~Cartan \cite{Cartan}  the 
notations
$P=-a_4,$ $ Q=-a_3,$ $R=-a_2,$ $S=-a_1,$ $ A=-L_1,$ $ B=-L_2$ were adopted, where
 $L_1$ and $L_2$ are the components of the projective curvature tensor.
In the work of R.Liouville \cite{Liouville} there were provided the semiinvariants $\nu_5,$ $w_1$, $i_2$, and the parameter $R_1$ (see review   \cite{BordagBandle}).  Relations between them and the pseudoinvariants $F$, $\Omega$, $N$ and  the component $H$ are as follows,
 $$
 F^5=\nu_5,\quad H=L_1(L_2)_x-L_2(L_1)_x+3R_1,\quad \Omega=-w_1-\frac{\nu_5 a_4}{L_1^3}-4\frac{(L_1)_x R_1}{L_1^3},\quad N=\frac{i_2}{3}.
 $$
Other pseudovectorial fields  and pseudoinvariants
 appeared firstly in  \cite{Sharipov1, Sharipov2, Sharipov3}.

\section{Solution of  equivalence problem for some Painlev\'e equations}

\subsection{Painlev\'e I equation}
The equivalence problem for this equation was effectively solved in paper \cite{Kartak1}.

{\bf Theorem 2}. {\it Equation (\ref{eq}) is equivalent to Painlev\'e I equation
$$\tilde y''=6 \tilde y^2+ \tilde x
$$ under the point transformations (\ref{zam}) if and only if the following conditions hold:
$F=0$ in (\ref{F}),  $A\ne 0$ or $B\ne 0$ in (\ref{alpha}),
$\Omega=0$ in (\ref{Omega1}), (\ref{Omega2}),
 $N=0$ in (\ref{N}),
 $W=0$ in (\ref{W}),
 $V=0$ in (\ref{V}),
 $\Theta\ne 0$ in (\ref{theta}),
 $L_1\ne 0$ in (\ref{L1}).
Invariants $I_1$ and $I_2$ in (\ref{inv7}) are functionally independent. The  point transformation is
$
\tilde x=  1/{\sqrt[5]{12 I_1}}, $ $\tilde y=\pm{\sqrt{I_2}}({\sqrt[5]{12^3}\sqrt[10]{I_1}}).$
}

Here the pseudoinvariant $W$ is introduced by (\ref{L1}), (\ref{theta}) and (\ref{phi1}), (\ref{phi2}),
\begin{equation}\label{W}
W= \nabla_{\boldsymbol \theta}L_1=(L_1)_{1.0}\theta^1+ (L_1)_{0.1}\theta^2-5L_1(\varphi_1\theta^1+\varphi_2\theta^2),
\end{equation}
and the pseudoinvariant $V$ is introduced by (\ref{L1}), (\ref{alpha}) and (\ref{phi1}), (\ref{phi2}),
\begin{equation}\label{V}
V=\nabla_{\boldsymbol \alpha}L_1= (L_1)_{1.0}B- (L_1)_{0.1}A-5L_1(B\varphi_1-A\varphi_2).
\end{equation}

{\bf Example.}  The equation
$$
\aligned
y''=&-\sin^3 y(6x\cos^2 y+\sin y)+
\frac 1x (-18x^3\cos^3 y\sin^2 y-3x^2\sin^3 y\cos y-2)y^\prime-\\
&-(18x^3\cos^4 y\sin y+3x^2\sin^2 y\cos^2 y)y^{\prime 2}
 -(6x^4\cos^5 y+x^3\sin y\cos ^3y+x)y^{\prime 3}
 \endaligned
 $$
 is equivalent to Painlev\'e I equation.
 The corresponding invariants and the  change of variables are
 $$
I_1=\frac 1{12}\frac{1}{x^5\sin^5 y}, \quad I_2=\frac{12x\cos^2 y}{\sin y},
\quad \tilde y=x\cos y, \quad \tilde x=x\sin y.
$$

{\bf Example.} The equation
$$
y''=6y^2+f(x)
$$
is equivalent to Painlev\'e I equation if and only if  $f(x)=mx+n,$ where $m,\, n$ are the constants, $m\ne 0$.
Let us check the conditions of Theorem 2:
$$
A=12,\quad B=0,\quad F=0,\quad \Omega=0,\quad N=0,\quad
W=\frac {f''(x)}{248832},\quad V=0,
$$
$$
\Theta=-\frac y{12},\quad L_1=-\frac{f'(x)}{20736},\quad I_1=\frac{f^{\prime 4}(x)}{12f^5(x)},\quad I_2=\frac{12y^2}{f(x)}.
$$

{\bf Examples.} Only these equations in the handbook by E.~Kamke \cite{Kamke} 
 $$
\begin{aligned}
&6.3 \qquad y''=6y^2+x,  \quad\text{Painleve I} \\
&6.5 \qquad y''=-ay^2-bx-c,\quad\text{with}\quad a,\, b\ne 0   \\
\end{aligned}
$$
are equivalent to  Painleve I equation.

\subsection{Painlev\'e II equation}
The equivalence problem for this equation was effectively solved in papers \cite{Kartak1} and \cite{Kartak2}.

{\bf Theorem 3. } {\it Equation (\ref{eq}) is equivalent to Painlev\'e II equation $$\tilde y''=2\tilde y^3+\tilde x\tilde y+a$$ with the parameter  $a=\pm J$ with $J=(4+10I_6-60I_3)/(50\sqrt{I_9})$ if and only if the following conditions hold: $F=0$ in (\ref{F}),  $A\ne 0$ or $B\ne 0$ in (\ref{alpha}),
$\Omega=0$ in (\ref{Omega1}), (\ref{Omega2}),
$M\ne 0$ in (\ref{M1}), (\ref{M2}), $I_1=18/5$ in (\ref{inv1}), $I_9\ne 0$ in (\ref{inv1-1}), invariant $J$ is constant. Among the invariants $I_3$, $I_6$, and $I_9$ from (\ref{inv1}), (\ref{inv1-1},  $k=3$) one can find two functionally independent. The  point transformation is
$\tilde y= {1}/{\sqrt[6]{2500 I_9}},$ $
\tilde x={5I_6}/{\sqrt[6]{2500I_9}}-3 J\sqrt[6]{2500 I_9}/2.
$}


{\bf Example.} The equation
$$
y''=(-2x^3-xy+a)y^{\prime 3}
$$
is equivalent to Painleve II equation with $\tilde a= a$. 
Let us check the conditions of Theorem 3,
$$
A=0,\; B=-12x,\, F=0,\,\Omega=0,\, M=\frac{288}{5},\, I_1=\frac{18}{5},\, J=\pm a,
$$
$$
I_3=\frac{2x^3+xy-a}{30x^3},\; I_6=\frac{2xy-3a}{10x^3},\; I_9=\frac{1}{2500x^6},\quad
\tilde y=x,\, \tilde x=y.
$$

{\bf Example.} The equation
$$
y''=y^3+f(x)y+g(x)
$$
is equivalent to Painleve II equation if and only if
$g(x)=c=const,$ $f(x)=mx+n,$ where $m,\, n$ are the constants, $m\ne 0$. Let us check the conditions of Theorem 3,
$$
A=6y,\; B=0,\, F=0,\,\Omega=0,\, M=\frac{72}{5},\, I_1=\frac{18}{5},\, J=\frac{g(x)y}{\sqrt{2}(f'(x)y+g'(x))},
$$
$$
I_3=\frac{y^3+f(x)y+g(x)}{15y^3},\quad I_6=\frac{2f(x)y+3g(x)}{5y^3},\quad
I_9=\frac{2(f'(x)y+g'(x))^2}{625y^8}.
$$

{\bf Examples.} Only these equations in the handbook by E.~Kamke \cite{Kamke} 
$$
\begin{aligned}
&6.6 \qquad y''=2y^3+xy+a, \quad\text{Painleve II,} \\
&6.8 \qquad y''=2a^2y^3-2abxy+b,\quad a,\, b \ne 0,\quad\text{where}\quad \tilde a=\frac 12,  \\
&6.9 \qquad y''=-ay^3-bxy-cy-d,\quad a,\, b \ne 0,\quad\text{where}\quad \tilde a=\frac{d\sqrt a}{b\sqrt{-2}}, \\
&6.142 \qquad y''=\frac{y'^2}{2y}+4y^2+2xy,\quad\text{where}\quad \tilde a=0, \\
&6.145 \qquad y''=\frac{y'^2}{2y}-\frac{ay^2}{2}-\frac{bxy}{2}, \quad a,\, b \ne 0,\quad\text{where}\quad \tilde a=0, \\
& 6.27\qquad y''=-ay'-bx^my^n,\quad n=3,\,\; m,\,b\ne 0,\quad\text{where}\quad \tilde a=0.\\
\end{aligned}
$$
are equivalent to  Painleve II equation.

\subsection{Painlev\'e III equation with 3 zero parameters}

A general form of the Painlev\'e equations III reads as
$$
y^{\prime\prime}=\frac     1y(y^{\prime})^2-\frac     1x
y^{\prime}+\frac 1x(ay^2+b)+cy^3+\frac dy.
$$
It is a 4-parameter family of equations, which we denote by $PIII (a, b, c, d)$.
If three 
of four  parameters vanish, all these equations are equivalent one to another. Referring to  work \cite{Hietarinta}, we write the change of variables:  1), 3): $x=\tilde x,\,y=1/{\tilde y}$, 2): $x=\tilde x^2/2,\,y=\tilde y^2$.
$$
PIII(0,b,0,0)\stackrel{1)}{\rightarrow}PIII(-b,0,0,0)\stackrel{2)}{\rightarrow}PIII(0,0,-b,0)\stackrel{3)}{\rightarrow}PIII(0,0,0,b),
$$
The equivalence problem for this equation was effectively solved in paper  \cite{Kartak2}.

{\bf Theorem 4. } {\it Equation (\ref{eq}) is equivalent to Painlev\'e III equation with 3 zero parameters  if and only if the following conditions hold: $F=0$ in (\ref{F}),  $A\ne 0$ or $B\ne 0$ in (\ref{alpha}),   $\Omega=0$ in (\ref{Omega1}), (\ref{Omega2}),
$M\ne 0$ in (\ref{M1}), (\ref{M2}), $I_1=3/5$,  $I_3=1/15$ from (\ref{inv1}).}

{\bf Example.} The equation 6.75 in the handbook by E.~Kamke \cite{Kamke} 
$$
y''=-\frac{2}{x}y'-e^y
$$
is not equivalent to Painleve III equation with 3 zero parameters. Let us check the conditions of Theorem 4,
$$
A=-e^y,\quad B=0,\quad F=0,\quad \Omega=0,\quad M=\frac{e^{2y}}{15},\quad
I_1=\frac 35,\quad I_3=\frac{1}{15}-\frac{4}{15x^2e^y}.
$$
{\bf Example.} The equation
$$
y''=f(x)y'-e^y
$$
is equivalent to Painleve III equation with 3 zero parameters if and only if function $f(x)$ is the solution of equation
$f^2(x)-f'(x)=0,$ hence $f(x)=1/(c-x),$ where $c$ is the constant. Let us check the conditions of Theorem 4,
$$
A=-e^y,\quad B=0,\quad F=0,\quad \Omega=0,\quad M=\frac{e^{2y}}{15},\quad
I_1=\frac 35,\quad I_3=\frac{1}{15}-\frac{2(f^2(x)-f'(x))}{15e^y}.
$$

{\bf Examples.} Only these equations in the handbook by E.~Kamke \cite{Kamke} 
$$
\aligned
&6.14\qquad\qquad y''=e^y,\\
&6.28\qquad\qquad y''=-ay'-be^y+2a,\quad b\ne 0,\qquad a=0\quad \text{or}\quad  a=-1,\\
&6.76\qquad\qquad y''=-\frac ax y'-be^y,\quad b\ne 0,\qquad a=0\quad \text{or}\quad a=1,\\
&6.77\qquad\qquad y''=\frac ax y'-bx^{4-2a}e^y,\quad b\ne 0,\quad a=1,\\
&6.83\qquad\qquad y''=-\frac{a(e^y-1)}{x^2},\quad a=-2,\\
&6.110(111)\quad y''=\frac{{y'}^2}y\pm\frac 1y,\\
&6.118\quad\qquad y''=\frac{{y'}^2}y-ay'-by^2+2ay,\quad b\ne 0,\qquad a=0\quad \text{or}\quad a=-1,\\
&6.127\quad\qquad y''=\frac{{y'}^2}y-by^2,\quad b\ne 0,\\
&6.172\quad\qquad y''=\frac{{y'}^2}y-a\frac{y'}x-by^2,\quad b\ne 0,\qquad a=0\quad \text{or}\quad a=1\\
\endaligned
$$
are equivalent to to Painleve III equation with 3 zero parameters.

\subsection{Equations PIII(0,0,0,0), PV(0,0,0,0), PVI(0,0,0,1/2)}

$$
{\small
\aligned
& PIII(0,0,0,0):\; y''=\frac 1y {y'}^2-\frac 1x y',\quad
PV(0,0,0,0):\;   y''= \left(\frac{1}{2y}+\frac{1}{y-1}\right){y'}^2-\frac 1x y',\\
&
 PVI(0,0,0,\frac 12): \; y''= \left(\frac   1y+\frac   1{y-1}+
\frac 1{y-x}\right)\frac{{y^{\prime}}^2}2-
 \left(\frac        1x+\frac
1{x-1}+\frac
1{y-x}\right)y^{\prime}+\frac{y(y-1)}{2x(x-1)(y-x)}.
 \endaligned}
$$
For the equations $PIII(0,0,0,0)$, $PV(0,0,0,0)$, $PVI(0,0,0,\frac 12)$ the conditions $A=0$ and $B=0$ hold where $A$ and $B$ are from (\ref{alpha})  then according to  Theorem 1 (Lie) they are equivalent to $y''=0.$

\subsection{Algorithm for Painlev\'e III-VI equations}

Equations Painlev\'e III-VI are  4-parameter families, they depend on parameters $a$, $b$, $c$, and $d$. As they are in the first case of intermediate degeneration, we should use formulas (\ref{inv1}), (\ref{inv1-1}) to calculate the invariants. It is a well-known fact that $I_2=0$ for all Painlev\'e equations, and therefore only two sequences generated by $I_1$ and $I_3$ are nontrivial.

The next useful fact is that all invariants are  rational functions of the variables ``$x$'' and ``$y$''.
Let us take the invariant $I_1$ (its formula is the simplest).
First we regard the symbol $I_1$ as a parameter in order %
to convert rational function into polynomial.
 This polynomial depend on variables ``$x$'' and ``$y$'', parameters $a$, $b$, $c$, $d$ (they are parameters of the equation), and the new parameter $I_1$. For example, invariant $I_1$ for the Painlev\'e III equation is
$$
I_1=\frac 35\frac {64(c^2y^8-22cy^4d+d^2)x^2-4(cy^6a+db+49(cy^4b+day^2))yx +a^2y^6-22ay^4b+b^2y^2}{(8cy^4x-8dx-yb+ay^3)^2}
$$
The assosiated polynomial is the following one,
$$
\aligned
P_1(x,y;a,b,c,d,I_1)=&5I_1(8cy^4x-8dx-yb+ay^3)^2 -(64(c^2y^8-22cy^4d+d^2)x^2-\\
&-4(cy^6a+db+49(cy^4b+day^2))yx+a^2y^6-22ay^4b+b^2y^2)=0
\endaligned
$$
We shall call ``$y$'' a ``higher'' variable.
In the same way we construct a polynomial $P_4(x,y;a,b,c,d,I_4)$ from the formula for invariant $I_4$. Using Buchberger algorithm (see \cite{Cox}), we reduce polynomials $P_1(x,y;a,b,c,d,I_1)$ and $P_4(x, y; a, b, c, d ,I_4)$ with respect to "higher" variable ``$y$''.
Aa a result we get a polynomial $Q_1(x;a, b, c, d, I_1, I_4)$ and a formula for variable $y=R_1(x;a, b, c, d, I_1, I_4)$, where $R_1$ is a rational function.
In the same way we construct a polynomial $P_7(x,y;a,b,c,d,I_7)$ by the invariant $I_7$.  Then we reduce it together with the polynomial  $P_1(x,y;a,b,c,d,I_1)$ and get a new polynomial $Q_2(x; a, b, c, d, I_1, I_7)$.

Now we reduce polynomials  $Q_1(x;a, b, c, d, I_1, I_4)$ and $Q_2(x;a, b, c, d, I_1, I_7)$ with respect to  the variable ``$x$''. We get a quantity $K(a, b, c, d, I_1, I_4, I_7)$ and a formula for variable $x=R_2(a, b, c, d, I_1, I_4, I_7)$, where $R_2$ is also a rational function.

Repeating this procedure as many times as necessary, we obtain a relation between the invariants $K(I_1, I_4,\dots)=0$ that is a {\it necessary condition} of the equivalence as well as the formulas for the parameters $a$, $b$, $c$, and $d$ and for the variables $x$ and $y$ via invariants. These formulas form the {\it sufficient conditions} and complete the solution.

The main difficulty of this method is a bulky form of these polynomials, and this is why at  present the equivalence problem is successfully solved only for Painlev\'e IV equation,  see \cite{Kartak3}. But the  final formulas are too complicated, so here we  present only the necessary conditions of the equivalence.

\subsubsection{Necessary conditions for Painlev\'e IV equation}
Equation Painlev\'e IV depends on two parameters $a$ and $b$,
$$
 PIV(a,b):\quad y^{\prime\prime}=\frac    1{2y}(y^{\prime})^2+\frac    32
y^3+4xy^2+2(x^2-a)y+\frac by.
$$
We introduce additional invariants using formulas (\ref{inv1}), (\ref{inv1-1}),
$$
J_1=\frac{5I_1}{72},\qquad  J_4=\frac{I_4}{2160},\qquad J_{10}=\frac{I_{10}}{12960},\qquad I_{10}=\frac{\nabla_{\alpha}I_4}{N}.
$$
{\bf Theorem 5.} {\it If equation (\ref{eq}) is equivalent to Painlev\'e IV equation under the transformations (\ref{zam}), then the following necessary conditions  hold:
$F=0$ in (\ref{F}),  $A\ne 0$ or $B\ne 0$ in (\ref{alpha}),  $\Omega=0$ in (\ref{Omega1}), (\ref{Omega2}), $M\ne 0$ in (\ref{M1}), (\ref{M2}), $Z\ne 0$ in (\ref{z}),
a) $K_0=0$ in (\ref{K0}) for PIV(a,0) equation;
b) $K_n=0$ in (\ref{Kn0}) for PIV(a,b) equation, $b\ne 0$.}
\begin{equation}\label{K0}
\aligned
K_0=& 4608J_1^4-3248J_1^3+808J_1^2+48000J_4J_1^2-\\
&- 16500J_4J_1-83J_1+1125J_4+125000J_4^2+3,
\endaligned
\end{equation}
and
\begin{equation}\label{Kn0}
\aligned
K=&2^{22}3^9J_1^9-2^{18}3^4 7229J_1^8+2^{14} 3^2(20412\cdot 10^3J_4+795377)J_1^7+\\
&+2^{10}\cdot 3\cdot 5(11664000J_{10}-293875200J_4-3170041)J_1^6+\\
&+2^9\cdot 3\cdot 5 (47628 \cdot 10^5J_4^2+347502500J_4-33816\cdot 10^3J_{11}+1574799)J_1^5+\\
&+2^8(550148750J_{10}+1701\cdot 10^7J_{10}J_4-15275925\cdot 10^4J_4^2-31879206254-\\
&-7217838)J_1^4+2^5(5312667+437746\cdot 10^4J_4+405 \cdot 10^7J_{10}^2+46305\cdot 10^8 J_4^3+\\
&+479194\cdot 10^6 J_4^2-1168733750J_{10}-129705\cdot 10^6 J_{10}4)J_1^3+2^2(-2157057-\\
&-337746700J_4+12948575\cdot 10^2J_{10}+6615\cdot 10^9J_{10}J_4^2+33184\cdot 10^7 J_{11}J_4-\\
&-697457\cdot 10^6 J_4^2-219765\cdot 10^8 J_4^3-23075\cdot 10^6 J_{10}^2)J_1^2+2^2\cdot 5 (9675\cdot 10^5 J_{10}^2+\\
&-17852625J_{10}+33823650J_4-847425\cdot 10^4J_{10}J_4-615125\cdot 10^6 J_{10}J_4^2+\\
&+8080625\cdot 10^5 J_4^3+11864525\cdot 10^3J_4^2+9261+7875\cdot 10^7 J_{10}^2J_4)J_1+\\
&+5^2(15435J_{10}-21609J_4-12027400J_4^2-16033\cdot 10^5J_4^3+11606\cdot 10^3J_{10}J_4+\\
&+5\cdot 10^7 J_{10}^3+343\cdot 10^8J_4^4+20875\cdot 10^5 J_{10}J_4^2-58\cdot 10^7J_{10}^2J_4-175\cdot 10^4J_{10}^2).
\endaligned
\end{equation}

{\bf Example.} The equation No. 34 from the book \cite{Ince} named "Painleve 34" equation
$$
XXXIV.\qquad y''=\frac{{y'}^2}{2y}+4 a y^2-xy-\frac{1}{2y},\quad a=const\ne 0
$$
is not equivalent to Painleve IV equation, although $K_n=0$ from (\ref{Kn0}).
Let us check the conditions of Theorem 5,
$$
A=6a-\frac{3}{2y^3},\quad B=0,\quad F=0,\quad \Omega=0,\quad M=\frac{9a(35+4ay^3)}{10y^5},\quad Z=0,\quad K_n=0.
$$

\section{Cases of Painlev\'e equation with non-trivial algebra of point symmetries}
In a general situation Painlev\'e equations have the trivial algebra of point symmetries. But in  some special cases the  dimension of the point symmetries algebra is  8, 2, or 1.

\begin{enumerate}
\item For equations $PIII(0,0,0,0)$, $PV(0,0,0,0)$, $PVI(0,0,0,\frac 12)$ we have $dim(Z)=8$.

\item For equation $PIII$ with 3 zero parameters we have $dim(Z)=2$, and  the operators are
$$
X_1=x\frac{\partial}{\partial x}+y\frac{\partial}{\partial y},\qquad X_2=
(x \ln x-2x)\frac{\partial}{\partial x}+y\ln x\frac{\partial}{\partial y}.
$$

\item For equation $PIII(0,b,0,d)$ or  $PIII(-b,0,-d,0)$ 
$$
\aligned
&PIII(0,b,0,d):\quad y''=\frac{y^{\prime 2}}{y}-\frac{y^{\prime }}{x}+\frac bx+\frac dy,\qquad X=x\frac{\partial}{\partial x}+y\frac{\partial}{\partial y},\\
&PV(a,b,0,0):\qquad y''=\left(\frac 1{2y}+\frac 1{y-1}\right)y^{\prime 2}-\frac {y'}x+
\frac{(y-1)^2}{x^2}\left(ay+\frac by\right),\quad X=x\frac{\partial}{\partial x}, 
\endaligned
$$
(they are equivalent under the transformations $x=\tilde x$, $y=1/\tilde y$, see \cite{Hietarinta}) and equation PV(a,b,0,0) we have $dim(Z)=1$.
\end{enumerate}

\section*{Acknowledgments}
I am grateful to Dr.
\ Ruslan Sharipov for useful remarks.

The work is partially supported by the Government of Russian Federation through Resolution No. 220, Agreement No. 11.G34.31.0042 and partially supported by the Russian Education and Science Ministry, Agreement No. 14.B37.21.0358.

\end{document}